\documentclass[12pt]{amsart}
\usepackage{ url, curves}

\newtheorem{theorem}{\noindent Theorem}

\newtheorem{definition}{\noindent Definition}

\urldef\vershik\url{vershik@pdmi.ras.ru}
\title[Kantorovich metric]{KANTOROVICH METRIC: INITIAL HISTORY AND LITTLE-KNOWN APPLICATIONS}
\begin{document}
\author{A.Vershik}
\address{Mathematical Institute of Russian Ac.Sci. St.Petersburg branch, 
Fontanka 27, St.Petersburg, 191023, Russia.
}
\email{vershik@pdmi.ras.ru}

\thanks{Grant RFBR 02-01-00093, NSh.2251.2003.1}

\begin{abstract}
{We recall the history of the transportation
(Kantorovich) metric and the Monge--Kantorovich problem.
We also describe several little-known applications: the first one
concerns the theory of decreasing sequences of partitions
(tower of measures and iterated metric), the second one
relates to Ornstein's theory of Bernoulli automorphisms 
($\bar d$-metric), and the third one is the formulation
of the strong Monge--Kantorovich problem in terms of 
matrix distributions.
Bibliography: $30$ titles. 
}
\end{abstract}

\maketitle



\section{ Introduction: the first papers on the trans\-por\-ta\-ti\-on problem}

The studies on the transportation problem could be called
a true pearl in the extremely rich scientific legacy of L.~V.~Kantorovich.
The beauty and naturalness of the formulation, the fundamental character of the
main theorem (optimality criterion), and, finally, the wealth of 
applications (some of them are realized, but new
applications keep on arising in areas that appear only now) -- 
all this allows us to 
place these studies among the classic mathematical works of the 20th century.
Undoubtedly,
the same words can be applied to the whole series
of papers on linear programming (from which the transportation problem
cannot be separated), which became the starting point for further
studies on mathematical economics, but here we will only dwell on the
remarkable role of what was later called the ``Monge--Kantorovich
problem'' and ``trans\-por\-ta\-ti\-on metric.''
\footnote{This metric has a dozen
of names known (one most used Vasserstein metric), because it has been rediscovered more than once and still keeps being
rediscovered. For many years I had to explain that many
metrics known in measure theory, ergodic theory, functional analysis,
statistics, etc., introduced in the 50s--80s, 
are special cases of the general definition of Kantorovich's transportation
metric. Many papers and books have appeared since then
(see, for example, [18]), 
but maybe it is only now (2004) 
that we can say that the publicity of
the main facts discovered by L.V.\ and his co-authors
matches their importance.}
In this introduction we do not intend to give a survey of this huge
subject; we will mention only the very first papers of L.V. and
his co-authors.

Apparently, L.V.\ conceived the formulation of the transportation problem
soon after he defined the general model 
of the production planning problem,
i.e., in the late 30s (the booklet [8]).   
However, if we 
judge from the date of the first publication, the transportation problem was born in 1942, with
the publication of the note [10],   
which later became famous.
The year itself predetermined the long road this paper 
had to walk to become known to specialists. 
The paper contains an explicit for\-mu\-la\-tion of
the general continuous transportation problem
on a compact metric space, the dual problem, and the optimality criterion. 
Later, in the small note [11]    
published  in {\it Uspekhi Mat.~Nauk}, Kantorovich
established a relation 
to Monge's problem of 
excavations and embankments,     
i.e., to the transportation problem on the Euclidean plane.
Since then, 
the general Kan\-to\-ro\-vich problem is sometimes called
the {\it Monge--Kantorovich problem} ({\it MK-problem} for short).
The next paper \cite{6.}, joint with a pupil of L.V., 
M.~K.~Ga\-vu\-rin,
was addressed rather to applied mathematicians and economists;
it contained a development of the method of potentials (a version of the method
of resolving multipliers
suggested by L.V.\ in 1939) for solving  the
finite-dimensional transportation problem.
Written long before publication, it ap\-pea\-red only in 1949,
and this delay was caused not by the wartime conditions,
but by the Soviet practice of that time, when
each scientific paper that even slightly touched economic
(not to mention socio-economic) problems had to go through
long and absurd censorship; besides, the paper was published not
in a journal, but in a special hard-to-reach volume.

Till 1956, i.e., during 18 years of existence of the new
mathematical economic theory, L.V.\ and his co-authors published
less than 10 papers on this subject (I remember  
G.~Sh.~Rubinshtein making up,
at my request, the complete
list of these papers in autumn 1956). Surely, not because
texts dealing with these problems were not written. L.V.\ had already
prepared a whole book on economics, whose destiny is an exact
and gloomy illustration of the system's attitude to
scientific studies that do not keep within obligatory schemes, rigid and
hence fruitless. A revised version of the book was not
published till almost twenty years later ([13]).    

In 1955--56, L.V.\ decided to ``open'' this topic;
he began to give public and special lectures, to popularize his theory.
The moment was chosen quite well. However, 
the wide distribution and acknowledgment of these studies
 were still a long way off.
One can read about all these events in the book [16]    
(in particular, in my paper [26]),      
but a detailed account of the whole story is still to be written.

Let us return to the transportation problem. The third 
important paper on this subject was the paper [15]   
by L.V.\ and his pupil and co-author G.~Sh.~Ru\-bin\-shtein.
It is this paper that contained an explicit definition of 
the norm in the space of measures related to the transportation
metric. The main observation was that the 
conjugate space to the space of measures with this norm is the space
of Lipschitz functions, and the optimality criterion is nothing more than
the dual definition of the norm as a supremum over the sphere
of the conjugate space. Before this paper it was not known whether
the space of Lipschitz functions is conjugate to any Banach space.
At that time (1956--57) I was interested in mathematical economics and
maintained close contacts with L.V.\
and G.~Sh.~Rubinshtein, and G.~Sh.\
described me in detail the stages of their work; in particular,
he said that L.V.\ was very satisfied by this interpretation
of the transportation problem. After this paper, the metric
is often called the {\it Kantorovich--Rubinshtein metric}.

Here it is worthwhile to make two remarks. Of course, the idea of
duality was contained from the very beginning both in the 
booklet of 1939 (the method of resolving multipliers) and in the note 
by L.V.\ in {\it Doklady Akad. Nauk SSSR} [9]   
-- the first paper devoted to comprehending
the relations between functional analysis and nonclassical
linear extremal problems (calculation of norms and extrema); 
it is worth noting that this was one more example
showing the utility of functional analysis for applications; see the 
paper [14],     
devoted to applications of linear programming
to computational mathematics, and the 
classical work by L.V.\ on the Newton method   [12].     
On the other hand, the technique that consists in
taking the objective function
as the norm in the space of right-hand sides of an extremal problem
(exactly this was suggested in [15])   
can be successfully applied to
many extremal problems (see, for example, [19, 28]).    
It was noted more than once that both classics of mathematical economics
of the 20th century -- von Neumann and Kantorovich --
came from functional analysis.

We cannot but mention that eventually, in course of 
development of the theory of nonclassical extremal problems, other
relations became obvious: to the theory of linear inequalities
and separability theory, 
Chebyshev ap\-pro\-xi\-ma\-tions and Krein's $L$-moment problem,
Weyl's studies on convex polytopes and convex geometry as a whole,
Bourbaki's theory of polars and com\-bi\-na\-to\-rics, etc.
\footnote{The lecture
course ``Extremal problems,'' which I taught for many years at the Department
of Mathematics and Mechanics of the Leningrad State University,
was compiled taking into account all these relations; in fact, it was a 
synthesis of functional analysis and the theory of extremal problems.
The textbook based on this course was not 
finished, but part of material was included in the 
textbook [1]     
written by my pupil A.~I.~Barvinok.}
Today we would include in this list ``tropical'' mathematics,
or max-plus algebra and impetuious developement of the applications
to differential equations, in particualr to Monge-Ampere equation,
hydrodynamics and so on (see references). We will not discuss here
those illuminated applications.

\section{ Basic definitions}

The transportation problem has been 
always holding a prominent position
among all problems of linear programming
due to its
general formulation and methods of solution. In what follows, I would like
to present several little-known ap\-pli\-ca\-ti\-ons
of the trans\-por\-ta\-ti\-on metric; but first let us recall the formulation
of the transportation problem.

\begin{definition}
{Let $(X,r)$ be a compact metric space, and let
$\mu_1$ and $\mu_2$ be two probability Borel measures on $X$. Consider the
{\it Monge--Kantorovich variational problem}  ({\it MK-problem} for short):
set
$$
k_r(\mu_1, \mu_2) =\inf_L \int r(x_1,x_2)\,dL,
$$ 
where $L$ runs over all Borel measures on $X \times X$
with marginal measures
$\mu_1$ and $\mu_2$.

The quantity $k_r(\mu_1, \mu_2)$ determines a metric on
 the simplex $V(X)$ of all probability measures on the compact space $X$;
it is called the {\it Kantorovich} (or {\it transportation}) {\it metric}
([10]).}    
\end{definition} 

\smallskip
\noindent{\bf Remark.} 
The measure $L$ is a ``plan of transportation''
of the distribution $\mu_1$ to the distribution $\mu_2$; the
integral means the cost of a given transportation plan,  and the infimum
(the Kantorovich metric) is achieved at the optimal plan.
\smallskip

\begin{theorem} {\rm (Kantorovich--Rubinshtein [15])}   
{\rm(1)} Consider the vector space $V_0(X)$ of all 
(not necessarily positive) Borel
measures $\nu$ with zero charge and finite variation (i.e.,
the positive part,
$\nu^{+}$, and the negative part, $\nu_{-}$, of $\nu$ have the
same finite variation) and define the Kantorovich--Rubinshtein norm
$||\nu||_k$ of an element $\nu \in V_0(X)$ as the Kan\-to\-ro\-vich
distance between the positive and negative parts
of $\nu$:
$$
||\nu||_k=k_r(\nu_{+},\nu_{-}).
$$
  Then the space of Lipschitz (up to additive constant) functions
  with the Lipschitz norm is the conjugate normed space to the space
   $V_0(X)$ with the norm $||.||_k$.

{\rm (2)} A plan $L$ in {\rm (1)} is optimal if and only if there exists
a Lipschitz function $U$ with Lipschitz constant $1$ such that
$$
U(x)-U(y)=r(x,y)
$$ 
almost everywhere with respect to the plan $L$.
\end{theorem}

We will omit the index $r$ in the notation 
$k_r$ if the metric $r$ is fixed, as well as the index $k$ in 
the notation $||.||_k$.

\smallskip\noindent{\bf Remark 1.}
The Kantorovich metric induces the weak topology on the simplex
of probability measures on the compact space $X$ ([15]).   

\smallskip\noindent{\bf Remark 2.}
In the framework of solution of the  finite-dimensional
trans\-por\-ta\-ti\-on problem, the optimal Lipschitz function $U$
is nothing more than the Kan\-to\-ro\-vich--Gavurin potential from [6].     
\smallskip

There is a huge number of difficult problems related to explicit cal\-cu\-la\-tion of
the Kantorovich metric for a given compact
space. For $\mathbb R^2$, this is the
classical Monge's problem on transportation of sand.
For $\mathbb R^1$, there is a good answer:
let $\nu_1$ and $\nu_2$ be two probability measures on $[0,1]$,
and let $r$ be the ordinary (Euclidean) metric; then
$k_r(\nu_1,\nu_2)=\int_0^1|\nu_1([0,t])-\nu_2([0,t])|\,dt$, i.e., the
Kantorovich metric is just the $L^1$-metric for distribution functions.
Apparently, there are no explicit formulas for $\mathbb R^n$, $n\ge2$.
Many papers are devoted to this problem; we will
 mention only the recent surveys [29, 30, 5].    

However, it makes sense to mention an essential idea, which has
ap\-peared recently and which plays a very important
role in modern applications to hyd\-ro\-dy\-na\-mics,
differential equations, and other areas (see [4,5,30] and references there);     
I mean the $p$-Kantorovich norms
(see [2]).     
Namely, the original definition of 
the Kantorovich metric (and Kantorovich norm)
resembles the definition of the $L^1$-norm; but we can also define
an analog of the $L_p$-norm
$$
k_p(\nu_1, \nu_2) =\inf_L\left[\int r(x_1,x_2)^p \,dL\right]^{1/p},
$$
where the infimum is taken, as before, over all transportation plans $L$ 
for a pair of probability
measures $(\nu_1, \nu_2)$, and the corresponding norm
$$
||\nu||_p = k_p(\nu_{+}, \nu_{-})
$$
for all $p \geq 1$. Of course, the original Kantorovich metric (the case $p=1$)
has more physical significance, but the case $p=2$ is much more
convenient from the technical and geometric point
of view. The corresponding variational
problem and Euler equation are simpler than in the case $p=1$, and the
results of [2]    
show that for a certain geometric transportation
problem, the Euler equation is the well-known Monge--Amp\`ere equation
(which {\it a priori} has nothing to do with the Monge--Kantorovich 
problem).

Let us mention another important special case, which is
sometimes also called the MK-problem; we will call it the {\it strong
MK-problem}. Namely, with the above notation, it is formulated as follows: {\it
to find
$$
{\bar k}(\mu_1,\mu_2)\equiv \inf_T \int r(x,Tx)\,d\mu_1(x),
$$
where the infimum is taken
over all measurable mappings $T$ such that $T\mu_1=\mu_2$.}

The existence of minimum in (2) is a very subtle
question. Of course, $\bar k(\mu_1,\mu_2) \geq k(\mu_1,\mu_2)$, and
the question of when the inequality becomes an equality is difficult and very
important. In the last section, we will present a new approach to
both problems.

 Among a huge number of applications of the Kantorovich
metric, I would like to mention only three examples, which are
little known to specialists in applications of this metric, 
yet are very important in
dynamical systems and functional analysis.

\section{The iterated Kantorovich metric and the tower of measures}

 We will begin with the notion of tower of measures, which was defined
in [20]   
and considered in more detail in [23], [39]   
Let $(X,r)$ be an
arbitrary compact metric space (say, the unit interval with the Euclidean metric).
We can consider a new compact space $V(X)$, the space of all
probability Borel measures on $X$, and supply it with the Kantorovich
metric. Thus we have defined a functor $F$ from the category of metric
compact spaces to itself: $F:X \mapsto V(X)$, $r\mapsto k_r$; it is
clear that $F$ sends each homeomorphism of a compact space $X_1$ to a
compact space $X_2$ to a homeo\-mor\-phism of $V(X_1)$ 
to $V(X_2)$.

Obviously, $(X,r)$ can be isometrically embedded into $(V(X),k_r)$
via the mapping $x \mapsto \delta_x$.

Let us iterate this procedure:
$$ 
(X,r) \longrightarrow (V(X),k_r) \longrightarrow
(V(V(X)),k_{k_r}) \longrightarrow \dots.
$$
Set $V^n=V(V^{n-1}(X))$ and $k_r^n=k_{k_r^{n-1}}$ and
introduce the notation $F_n$ for the mapping $(V^{n-1},k_r^{n-1})
\longrightarrow (V^n,k_r^n)$.

We can consider the {\it inductive limit} of this sequence
of metric spaces with isometric embeddings:
$$
(V^{\infty},k_r^{\infty})\equiv \mbox{indlim}_n((V^n,k_r^n), F_n).
$$
This inductive limit (a metric space)
is called the {\it infinite tower of measures}; it
plays a crucial role in the
theory of filtrations of $\sigma$-fields generated by random processes
and its various applications.

On the other hand, for $n\ge2$ there is a natural projection
$$
P_n:V^n \longrightarrow V^{n-1}, \qquad P_n(\mu)=\bar\mu, 
$$
where $\bar \mu$ is the barycenter of the measure $\mu$,
which is well defined
for measures on affine compact spaces (thus the projection is defined
for $V^n$, $n\geq 2$), and we have the sequence
$$
(V^1(X),k_r)\longleftarrow (V^2(X),k_r^2) \longleftarrow \dots. 
$$
Thus we obtain the {\it projective limit}
$$
{\overline V}^{\infty}\equiv \mbox{projlim}_n (V^n(X),P_n).
$$
Since $P_n F_n= I_{n-1}$,
the inductive limit $V^{\infty}$ is naturally embedded into the
 projective limit:
$$ 
V^{\infty} \subset {\overline V}^{\infty};
$$
but, in contrast to the case of inductive limit, on the projective limit there
is 
 no natural metric.\footnote{Inductive systems having projections
 that are the right inverses to the embeddings can be called 
 {\it indo-projective} systems; they appear quite often.}

 The main application of this tower of measures is as follows.
Assume that we have a ``metric triple''
 $(X,r,\mu)$, i.e., a measure space with a metric or semimetric,
and a decreasing sequence of measurable partitions of this space
(discrete filt\-ra\-ti\-on) 
 $\{\xi_n\}$, $n=0,1, \dots$; here $\xi_0$ is trivial and
 $\xi_n > \xi_{n+1}$.

First consider one partition $\xi$; for almost all points
 $a \in X/\xi$ of the quotient space with respect to this partition, there is a well-defined
 conditional measure on the element of  $\xi$
 corresponding to $a$. We regard it as a
 measure on $(X,r)$; thus we have a mapping $f_{\xi}:X/\xi \rightarrow V(X,r)$,
 which sends almost every point $a \in X/\xi$
 to a (conditional) measure on $(X,r)$.
 It is convenient to regard this mapping as a function from
 $(X,\mu)$ to $V(X)$.

 Now define a metric (or semimetric) on $X/\xi$ as follows:
 for almost all pairs of points $a,b \in X/\xi$, define
 the distance between them as the Kantorovich distance between
the corresponding  conditional measures.
 Thus we have defined a metric (or semimetric) on a subset
 of full measure in the quotient space $X/\xi$; it can also be regarded
 as a semimetric on the original space $(X,\mu)$.

 Apply this process to the decreasing sequence
 of partitions $\{\xi_n\}$: start from $\xi_1$,
 then define a metric on $X/\xi_1$, a mapping
 $f_1:X\rightarrow V(X,r)$, and a partition
 $\xi_2/\xi_1$; now we have a mapping from
 $X/\xi_2$ to $V^2(X)$, a new metric on $X/\xi_2$, and a
 map $f_2:\rightarrow V^2(X,r)$.

 Continuing this process, we obtain mappings $f_n$
 from $(X,\mu)$ to the iterated spaces $V^n(X,r)$, or
 to the inductive limit $(V^{\infty},k_r^{\infty})$.

 One of the main results of the theory of decreasing
 sequences ([20], [23]) 
is the following theorem.
 
 \begin{theorem}
A decreasing homogeneous sequence of measurable partitions
 is {\it standard} (see [20, 23]      
for definitions) if  and only if the
 sequence of measures $f_n*\mu$ (in other words, the sequence of the
 distributions of the mappings $f_n$ with respect to the measure $\mu$), 
 regarded as a sequence
 of measures on the inductive limit $(V^{\infty},k_r^{\infty})$,
 tends to a $\delta$-measure.
 \end{theorem}

 A discussion of these subjects can be found in [23]   
and in forthcoming papers.

\section{The Kantorovich metric in Ornstein's theory}

  In the early 70s, Donald Ornstein solved a long-standing problem
  in ergodic theory: he gave necessary and sufficient conditions
  on a discrete-time stationary random process under which the shift
  in the space of trajectories of this process is isomorphic to a Bernoulli
  shift; using this result, he proved that the Kolmogorov entropy 
  is a complete invariant
  of Bernoulli shifts ([17]).    
We will formulate the main theorem of Ornstein's theory 
in order to illustrate the role of the Kantorovich metric, which was
  rediscovered by Ornstein (he called it the $\bar d$-metric).
  
Assume that the state space $S$ of a stationary process is finite
  and $\mu$ is the stationary measure on $S^{\mathbb Z}$ generated by this process.
  The question is formulated as follows: when there exists an isomorphism (in
  the measure-theoretic sense) of the Bernoulli space ${S'}^{\mathbb Z}$ with
  product measure and the space $(S^{\mathbb Z}, \mu)$ that commutes with the shift.
  This is the well-known iso\-mor\-phism problem in ergodic theory.
  It is clear that the criterion of existence of such an isomorphism must
  be expressed in terms of the rate of
  decrease
  of the correlation between the past and the future of the process.
  There are many known conditions of this type, which are 
  sometimes called
  ``mixing conditions.'' Most of such conditions known in the theory
  of stationary processes are too strong (Kolmogorov's, Rozenblatt's, Ibragimov's
  conditions, etc.). It turned out that the right notion is related to the
  Kantorovich metric
  on the space of words with the Hamming metric -- this was
  discovered by D.~Ornstein. Our interpretation slightly
  differs from the original one, but is closer to the previous context
  (see [29]).

  Let $\{\xi_n\}$, $n \in\mathbb Z$, be a stationary random process with finite
  state space $S$ and shift-invariant measure $\mu$ on $S^{\mathbb Z}$.
Consider the ``past'' of the process:
  ${\mathcal P}=\prod_{-\infty}^0 S$; the projection of $\mu$ to $\mathcal P$ will be 
  denoted by $\mu^-$.
  Fix a point $x^-=(x_0,x_{-1},x_{-2}, \dots) \in \mathcal P$
  and consider the conditional distribution on the
  $n$-future given a  fixed past $x^-$:
$$
P_n(x_1,x_2,\dots, x_n|x^-);
$$
  this is a measure on the $n$-future $S^n$ defined for
  almost all points $x^- \in {\mathcal P}$;  it is an element
  of $V(S^n)$, thus we have a mapping $F_n: {\mathcal P} \to V(S^n)$
  defined almost everywhere.

 Consider the Hamming metric on $S^n$: 
 $$
 h_n(x,y)=
 \frac1n\#\{i \in (1, \dots, n):x_i\neq y_i\},
 $$
 where $x=(x_1, \dots, x_n),
 y=(y_1,\dots, y_n)\in S^n$ and $\#$ stands for the 
 number of points in a set;
 and let $k_{h_n}$ be the Kantorovich metric  on the space 
 $V(S^n,h_n)$ of measures on the $n$-future.

 \begin{theorem} {\rm[17, 24]}
Consider a stationary process $\{\xi_n\}$, $n \in \mathbb Z$,
 and the right shift in the space of realizations
 generated by this process.
An invertible encoding of this shift into a Bernoulli
 shift (in other words, a measure-preserving isomorphism of the
 shift in the space of realizations of the process and a
 Bernoulli shift) exists if and only if
 $$
\lim_{n\to \infty} \iint\limits_{x^-\in {\mathcal P},\, y^- \in
 {\mathcal P}} k_{h_n}(P(*|x^-),P(*|y^-))\;d\mu^-(x^-)d\mu^-(y^-)=0
$$
 (the integral of the value of the
 Kantorovich metric for the pair of conditional measures
 corresponding to a pair of points from ${\mathcal P}\times{\mathcal P}$
 with respect to the product measure $\mu^-\times\mu^-$).
 \end{theorem}
 
 The literal meaning of the above condition is very transparent: it
 means that the conditional distribution on the future given a
 fixed past asymp\-to\-tical\-ly does not depend on the past; roughly
 speaking, there is only one type of distribution on the future; but
 a more precise sense of these words essentially depends on the choice of
 a metric on the space of realizations of the process (we should take the
 Hamming metric) and a metric on the spaces of measures (here we should
 use the Kantorovich metric); in general, the conclusion of the theorem will be
 false if we replace the Kantorovich metric by some other one (for
 example, by the variation metric).

 The last formulation also motivates the definition of the so-called
 {\it se\-con\-da\-ry entropy} of a stationary process (see [24]).   
 Define $M_n^+$ as the image of the measure $\mu^-$ (see above) under the
 mapping $F_n:{\mathcal P}\to V(S^n,h_n)$;  this is a measure on
 $V(S^n,h_n)$. In the case of Bernoulli automorphisms,
 by Ornstein's theorem, the measure $M_n^+$
 tends to a $\delta$-measure as $n\to\infty$. But for a
 general Kolmogorov stationary process (K-automorphism), this is not
 the case. More precisely, if the automorphism is not a Bernoulli
 automorphism, then the limit exists, but is not a $\delta$-measure.
 Thus it is natural to introduce a characteristic of the limiting measure.
 Namely, we may consider the so-called $\varepsilon$-entropy of the measure
 $M_n^+$. This notion also uses the Kantorovich metric.  For 
 an arbitrary Borel probability
 measure $\nu$ on a metric space $(X,d)$, the $\varepsilon$-{\it entropy}
 $h_{\varepsilon}(\nu)$ (as a function of $\varepsilon$) is defined as follows:
 $$
h_{\epsilon}(\nu)= \inf\{H(l):k_d(l,\nu)<\epsilon \},
$$
 where the infimum is taken over all discrete measures $l$ on $(X,d)$ and
 $H(l)$ is the ordinary entropy of a discrete measure: $H(l)=-\sum l_i
 \log l_i$, $l=(l_1,\dots, l_n)$, $\sum_i l_i=1$, $l_i\geq 0$, $i=1, \dots, n$.

 The asymptotic of $h_{\varepsilon}(M_n^+)$ with respect to $n$ is called the 
 {\it secondary entropy}
 of the process. An open problem: what kind of asymptotic behavior
 can appear?
 Presumably, the secondary entropy is a metric invariant of K-automorphisms.

\section{Application to the classification of metric spaces}

  Consider a Polish (=metric, complete, separable) space with a
  Borel probability measure. We call such a space a {\it metric triple} (another
  term is an $mm$-space [7]).    
  Two triples $(X, \rho, \mu)$ and  $(X', \rho', \mu')$  are isomorphic
  if there exists a mapping $T:X \rightarrow X'$ that is an isometry and
  preserves the measures:  $\rho'(Tx, Ty)=\rho (x,y)$ and
  $T\mu=\mu'$.

  We regard the metric as a measurable function of two variables:
  $$
\rho:X \times X \longrightarrow \mathbf R. 
$$
 (The theorem below is true for an arbitrary symmetric
  measurable function $\rho$, not necessarily a metric.)

  Let $X^\infty$ be the product of infinitely many copies of the space $X$.
  Define a mapping
  $$
 F:X^\infty \longrightarrow M_{\infty}(\mathbb R)
$$
  from $X^\infty$ to the set of symmetric matrices
  as follows:
  $F(x,y)= \{r_{i,j}\}_{i,j=1}^{\infty}$,
  where $x=(x_1,x_2,\ldots)$
  and $r_{i,j}=\rho(x_i,x_j)$.

  Let us denote the image of the measure $\mu^{\infty}$
  under the mapping $F$ by  $F(\mu) \equiv D_{\rho}$;
  the measure $D_{\rho}$ on $M_{\infty}(\mathbf R)$
  will be called the {\it matrix distribution} of the function $\rho$.

  In [25],     
we considered and classified general (nonsymmetric)
  measurable functions $f(x,y)$ of two variables on the space
  $(X \times X, \mu \times \mu)$ up to mappings of the form $T_1 \times T_2$,
  where $T_1$ and $T_2$ are measure-preserving automorphisms of $(X,\mu)$.
  We also defined the notion of matrix distribution for this case; it
  is a complete invariant for so-called pure functions.
  But now we need another classification. We also consider
  arbitrary measurable
  (nonsymmetric) functions $f$  on the space $(X\times X, \mu \times\mu)$,
  where $(X,\mu)$ is a Lebesgue space with continuous measure,
  but we classify them up to mappings of the form $T \times T$, 
  where $T$ is an automorphism
  of $(X, \mu)$ (in other words, $T_1=T_2$). Namely, define a mapping
  $$
F_f:X^\infty \longrightarrow M_{\infty}(\mathbb R),
$$
  where $F_f(x)=\{f(x_i,x_j)\}_{i,j=1}^{\infty}$ and
  $x=(x_1,x_2, \dots)\in X^{\infty}$; here $M_{\infty}(\mathbb R)$ is the
set of
  arbitrary (not necessarily symmetric) matrices.
  The $F_f$-image of the measure  $\mu \times \mu$,
  which is a measure on $M_{\infty}(\mathbb R)$, is called the {\it
  symmetric matrix distribution} of the function $f$
  and denoted by $D^s_f$.

  \begin{theorem} {\rm (Gromov [7], Vershik [25])}

{\rm(1)} Two metric triples $(X, \rho, \mu)$ and  $(X', \rho', \mu')$
are isomorphic if and only if their matrix distributions coincide:
$$ 
D^s_{\rho}=D^s_{\rho'}.
$$ In other words, the matrix distribution   of the metric
  is a complete invariant of a metric triple.

  {\rm(2)} {\rm(}Vershik [25]{\rm)}.     
The symmetric matrix distribution $D^s_f$
  of a measurable function $f(\cdot,\cdot)$ of two variables
  is a complete metric invariant of the function regarded up to 
  automorphisms of
  the form $T \times T$, where $T$ is an automorphism of $(X,\mu)$.
  \end{theorem}

   Now we apply this classification to MK-problems.
   Let $X$ be a compact metric space with metric $\rho$;
   we want to ``transport'' a Borel probability measure $\mu_1$ to
   another Borel probability measure $\mu_2$.  Thus we have
   two metric triples: $(X,\rho,\mu_1)$ and $(X,\rho,\mu_2)$.
   It is more convenient to reduce the problem to a more
   symmetric form and to have one metric triple. Let us consider
   only continuous measures; then we can choose a
   measure-preserving isomorphism
   $S:(X,\mu_2) \to (X,\mu_1)$. Let $f(x,y)=\rho(x,Sy)$, 
   so that $f$ is a nonnegative measurable (in general, nonsymmetric)
   function of two variables -- the ``shifted metric.''
   Now we can consider only one measure $ \mu_1\equiv \mu $
   and the function $f$ on the space $(X \times X, \mu\times\mu)$.

   In terms of the shifted metric, the MK-problem
   can be formulated as follows:
{\it    to find
$$
k\equiv\inf_L \int f(x_1,x_2)\,dL,
$$
   where $L$ runs over all Borel measures on the product $X \times X$
   with both marginal measures equal to the
   measure $\mu$}; thus
   $L$ belongs to the set of {\it bistochastic measures}, or, in other
   words, $L$ is an element of the semigroup of {\it polymorphisms with invariant
   continuous measure} $\mu$ (see [22])  
for definitions). 
   Thus the MK-problem turns into 
   a variational problem on the convex set of bistochastic measures
   (or on the semigroup of polymorphisms).

   The strong MK-problem reads as follows: {\it
   to find
$$
\bar k \equiv\inf_T \int f(x,Tx)\,d\mu(x),
$$
   where $T$ runs over all $\mu$-preserving transformations of
   $(X,\mu)$}. In this case, we have a variational problem on the group of measure-preserving
   transformations.

   Now we can apply the above-defined
   symmetric matrix distribution $D^s_f$ 
   of the function $f$ regarded
   as a measurable function (shifted metric)
   on the space $(X \times X, \mu \times \mu)$.
   Since $D^s_f$ is a complete invariant of the triple $(X,f,\mu)$, {\it all
   properties of the (ordinary and strong) MK-problem 
   can be expressed in terms of $D^s_f$} as a measure on the space of matrices
   $M_{\infty}(\mathbb R)$. But this means that we have a {\it random matrix}
   with distribution $D^s_f$,  which we can use for analysis of the
   problem. Here we describe only one example of applying this approach.

   Let $r=\{r_{i,j}\}_{i,j=1}^{\infty}$ be a random matrix with 
   distribution $D^s_f$. The new version of the MK-problem reads as follows.
   Choose a random matrix $r$, for each $n$ consider the ordinary finite
   transportation problem, and define
$$
k_n(r) \equiv \inf_{l} \sum_{i,j=1}^n l_{i,j} r_{i,j},
$$
   where $l=\{l_{i,j}\}_{i,j=1}^n$  is a bistochastic matrix
   (i.e., $\sum_{i=1}^n l_{i,j}=\sum_{j=1}^n l_{i,j}=1$,
   $l_{i,j}\geq 0$ for all $i,j=1, \dots, n$) and
   $r_n=\{r_{i,j}\}_{i,j=1}^n$ is the $n$-fragment of $r$
   (the random matrix constructed from the shifted metric 
   as described above). Thus
   $k_n(r)$ is a random variable that depends on the random matrix $r$.

\begin{theorem} In the previous notation,
$$
\lim_{n \to \infty} k_n(r) = k\quad \text{in measure}\quad D^s_f,
$$ 
where $k$ is the solution of the original MK-problem,
    i.e., the sequence of random variables $k_n(r)$ converges
    {\it in measure} $D^s_f$ to the solution  of the MK-problem.
   \end{theorem}

   A natural conjecture:
   {\it for almost every choice of the matrix $r=\{r_{i,j}\}$ 
   with respect to the measure
   $D^s_f$, the same assertion is true:
$$
D^s_f\{r:\lim_{n \to \infty} k_n(r)=k\}=1,
$$
   which means that $k_n(r)$ converges to $k$ with probability one
   with respect to the choice of the matrix $r$ according to the measure $D^s_f$.}

   Note that we approximate the MK-problem with the simplest finite-dimensional
   problem of linear programming -- the allocation problem. 
   By the Birkhoff--von Neumann
   theorem, the solution of this problem 
   is a permutation, i.e., an element of the symmetric group, or an
   extreme point of the convex set of bistochastic matrices (the
   so-called Hungarian
   polytope). Nevertheless, the question of when the 
   strong MK-problem has a solution
   and how it can be approximated by permutations is more involved.

   The theorem and conjecture given above are typical
   for applications of our method to various problems with integral kernel:
   we obtain a probabilistic approximation of a functional or
   variational problem using a random
   choice of values of the function. We will return to this elsewhere.
\medskip

Partially supported by the RFBR, project 02-01-00093, and 
the President of Russian Federation grant for 
support of leading scientific schools NSh-2251.2003.1.
\medskip

   Translated by A.~M.~Vershik and N.~V.~Tsilevich.

\end{document}